\documentclass[12pt,leqno]{article}
\newcommand{\pbk}{\pagebreak}

\newtheorem{thm}{Theorem}[section]
\newtheorem{lem}[thm]{Lemma}
\newtheorem{pro}[thm]{Proposition}
\newtheorem{cor}[thm]{Corollary}
\newtheorem{Def}[thm]{Definition}
\newcommand{\bdfn}{\begin{Def} \rm}
\newcommand{\edfn}{\end{Def}}
\newcommand{\bthm}{\begin{theorem} \rm}
\newcommand{\ethm}{\end{theorem}}
\newcommand{\pf}{\noindent{\bf Proof. }}
\newcommand{\lef}{\left \langle}
\newcommand{\rig}{\right \rangle}
\newcommand{\DPP}{Dunford-Pettis property }

\title{\Large \bf An operator summability of sequences in Banach spaces}
\author{\normalsize {\bf Anil K. Karn$^{\dagger \star}$ and Deba P. Sinha$^{\ddagger }$}\\[-6pt]}

\date{}
\begin{document}

\maketitle
\footnote{
\begin{description} 
\item[$\dagger$] School of Mathematical Sciences, National Institute of Science Education and Research, Institute of Physics Campus, P.O. Sainik School, Bhubaneswar 751005, India \ {\it E-mail}: {\texttt anilkarn@niser.ac.in}
\item[$\ddagger$] (Lately at) Department of Mathematics, Dyal Singh College (University of Delhi),
Lodi Road, New Delhi 110003, India
\item [$\star$] Corresponding author.
\item [2010  Mathematics Subject Classification:] Primary 46B20; Secondary 46B28, 46B50.
\item [Keywords and phrases:] Weak and norm summable sequences, Compact and limited sets, compact and absolutely summing operators, the Dunford-Pettis property and the Gelfand-Phillips property.
\end{description}
}
\thispagestyle{empty}

\begin{abstract}
Let $1 \leq p <\infty$. A sequence $\lef x_n \rig$ in a Banach space $X$ is defined to be $p$-operator summable if for each $\lef f_n \rig \in l^{w^*}_p(X^*)$, we have $\lef \lef f_n(x_k)\rig _k \rig _n \in l^s_p(l_p)$. Every norm $p$-summable sequence in a Banach space is operator $p$-summable, while in its turn every operator $p$-summable sequence is weakly $p$-summable. An operator $T \in B(X, Y)$ is said to be $p$-limited if for every $\lef x_n \rig \in l_p^w(X)$, $\lef Tx_n \rig$ is operator $p$-summable. The set of all $p$-limited operators form a normed operator ideal. It is shown that every weakly $p$-summable sequence in $X$ is operator $p$-summable if and only if every operator $T \in B(X, l_p)$ is $p$-absolutely summing. On the other hand every operator $p$-summable sequence in $X$ is norm $p$-summable if and only if every $p$-limited operator in $B(l_{p'}, X)$ is absolutely $p$-summing. Moreover, this is the case if and only if $X$ is a subspace of $L_p(\mu )$ for some Borel measure $\mu$.
\end{abstract}

\pbk

\section {\large \bf Introduction}

Let $X$ be a Banach space, $\lef x_n \rig$ a sequence in $X$ and $1 \leq p < \infty$. We say that $\lef x_n \rig$ is (norm) $p$-summable in $X$ if
$\sum _{n=1}^{\infty}{\|x_n\|^p} <\infty$. If $\sum _{n=1}^{\infty}{|f(x_n)|^p}< \infty$, for all $f\in X^*$, then we say that $\lef x_n \rig$ is weakly $p$-summable in $X$. It is easy to note that a norm $p$-summable sequence is always a weakly $p$-summable, while the converse, in general, is not true. In fact in a Banach space $X$ every weakly $p$-summable sequence is norm $p$-summable if and only if $X$ is finite dimensional. These two types of summability were used by Grothendieck \cite{G1} to introduce the operator ideal of absolutely summing operators (for $p=1$), further generalized by Piestch \cite{P} who defined the operator ideal of absolutely $p$-summing operators for all $1 \leq p < \infty$. These operator ideals have been studied extensively in the literature.

Let $l^s_p(X)$ denote the set of all norm $p$-summable sequences and $l^w_p(X)$ that of all weakly $p$-summable sequences in $X$. Then these two sets become Banach spaces under suitable norms. More precisely, $l^s_p(X)$ can be identified as the 'countable $p$-direct sum' of $X$; similarly, $l^w_p(X)$ can be shown to be isometrically isomorphic to the space $B(l_{p'}, X)$ of operators if $p>1$ (here $p'$ is the harmonic conjugate of $p$, i.e. $\frac{1}{p} + \frac{1}{p'} = 1$) and $l_{p'}$ is replaced by $c_0$ when $p=1$.

In this paper we introduce a new kind of summability of sequences in Banach spaces using the notion of $p$-summing operators and call it the operator $p$-summability (definition below). This notion crops up naturally while extending the idea of limited sets to a $p$-level. In general, this type of summability of sequences is different from both weak and norm summability. In this paper, we investigate Banach spaces for which this type of summability coincides either with weak or with norm summability. For the first type of Banach spaces in question, we encounter a $p$-level of Dunford-Pettis property whereas for the other we are encouraged to introduce the notion of a $p$-level of Gelfand-Phillips property. The later type of Banach space ultimately reduces to subspaces of $L_p(\mu )$ for some Borel measure $\mu$.

Example of a Banach space can be constructed for which the operator $p$-summability is different from both norm as well as weak $p$-summabilities.

\section{An operator summability}

A non-empty subset $S$ of a Banach space $X$ is said to be limited if for every weak$^*$-null sequence $\lef f_n \rig$ in $X^*$ (i.e., $lim_{n\to \infty}f_n(x)=0$, for all $x\in X$), $f_n \to 0$ uniformly on $S$. Alternatively, given a weak$^*$-null sequence $\lef f_n \rig$ in $X^*$ there is an $\lef \alpha _n \rig \in c_0$ such that $|f_n(x)|\leq \alpha _n$ for all $x\in S$ and all $n\in {\bf N}$. We can extend this idea to the '$p$-sense' in the following way. We define a subset $S$ of $X$ to be $p$ limited in $X$ $(1 \leq p < \infty)$ if for every weak$^*$-$p$-summable sequence $\lef f_n \rig$ in $X^*$ (i.e., $\sum _{n=1}^{\infty}{|f_n(x)|^p< \infty}$ for all $x \in X$) there is an $\lef \alpha _n \rig \in l_p$ such that $|f_n(x)| \leq \alpha _n$ for all $x \in S$ and $n \in \bf{N}$.

The history of limited sets originated from the following error of Gelfand \cite{G} : A set $S$ in Banach space $X$ is compact if and only if every weak$^*$- null sequence in $X^*$ is uniformly null on $S$. Clearly, every compact set has this property. However, Phillips \cite{Ph} came out with an example of a non-compact set with the above property, i.e., of a limited non-compact set. The authors \cite{SK} recently, studied the concept of
$p$-compact sets for $1 \leq p < \infty$. It is interesting to note that the above mentioned analogy carries over to $p$-level too. First we show that $p$-compact sets are $p$-limited.

We begin with some definitions. For $x = \lef x_n \rig \in l_p^w(X)$, we define an operator $E_x : l_{p'} \to X$ given by $E_x(\alpha)= \sum _{n=1}^\infty {\alpha _n x_n}$, $\alpha =\lef \alpha _n \rig \in l_{p'}$. Then $E_x \in B(l_{p'}, X)$. Moreover, in this identification $l_p^w(X)$ is isometrically isomorphic to $B(l_{p'}, X)$. For $p=1$, $l_{p'}$ is replaced by $c_0$. We say that $K \subset X$ is relatively $p$-compact
if there is an $x= \lef x_n \rig \in l_p^s(X)$ such that $K \subset E_x (Ball(l_{p'}))$. Similarly, $K \subset X$ is said to be (relatively) weakly $p$-compact if there is an $x = \lef x_n \rig \in l_p^w(X)$ such that $K \subset E_x (Ball(l_{p'}))$.

\begin{pro}
Let $1 \leq p < \infty$ and $X$ a Banach space. Then every $p$-compact subset of $X$ is $p$-limited.
\end{pro}
\pf
Let $K \subset X$ be a $p$-compact and $\lef f_n \rig \in l_p^{w^*}(X^*)$. There is an $x = \lef x_k \rig \in l_p^s(X)$ such that $K \subset E_x(Ball(l_{p'}))$. Then
$$\sum _{k=1}^\infty {\sum _{n=1}^\infty {|f_n (x_k)|^p}} \leq (\|\lef f_n \rig \|_p^{w^*})^p \sum _{k=1}^\infty {\|x_k\|^p},$$
so that
$$\sum _{n=1}^\infty {\sum _{k=1}^\infty {|f_n (x_k)|^p}} = \sum _{k=1}^\infty {\sum _{n=1}^\infty {|f_n (x_k)|^p}} \leq (\|\lef f_n \rig \|_p^{w^*} \|x\|_p^s)^p.$$
Set $(\sum _{k=1}^\infty {{|f_n (x_k)|^p)}^{\frac{1}{p}}} = \alpha _n$ for all $n$ so that $\lef \alpha _n \rig \in l_p$. Now if $z \in K$, then $z = \sum _{k=1}^\infty{\beta _k x_k}$ for some $\lef \beta _k \rig \in Ball(l_{p'})$, and for each $n$ we have 
$$|f_n(z)|= |\sum _{k=1}^\infty {\beta _k f_n(x_k)}| \leq (\sum _{k=1}^\infty |\beta _k|^{p'})^{\frac {1}{p'}} (\sum _{k=1}^\infty {|f_n (x_k)|^p})^{\frac {1}{p}} \leq \alpha _n.$$
Hence, $K$ is $p$-limited. \hfill $\triangle$
\vskip 8pt plus 0fill
Next we observe certain facts about $p$-limited sets.

\begin{pro}
Let $A$ and $B$ be two subsets of a Banach space $X$. \par
{\leftskip=6mm
\item{(a)} If $B$ is $p$-limited and $A \subset B$, then $A$ is also $p$-limited.
\item{(b)} If $A$ is $p$-limited then $\overline{A}$ is $p$-limited.
\item{(c)} If if $A$ and $B$ are $p$-limited sets then so are $A \cup B$, $A + B$ and $A \cap B$.
\item{(d)} If $A$ is $p$- limited and $T \in B(X, Y)$, then $T(A)$ is $p$-limited in $Y$.\par}
\end{pro}

\pf
Suppose $A$ is $p$-limited. We prove (b). Let $\lef f_n \rig \in l_p^w{^*}(X^*)$. Then there is an $\lef \alpha _n\rig \in l_p$ such that $|f_n (x)|\leq \alpha _n$ for all $x \in A$ and $n \in \bf{N}$. Let $x \in \overline{A}$. Then there is an $\lef x_k \rig$ in $A$ such that $x_k \to x$. Thus for each $n$, $f_n (x_k) \rightarrow f_n (x)$. Fixing $n$ we have $|f_n (x_k)| \leq \alpha _n$ for all $k$. It follows that $|f_n (x)| \leq
\alpha _n$ for all $n$, so that $\overline A$ is $p$-limited. Thus (b) follows. The proofs of (a), (c) and (d) are immediate. \hfill $\triangle$

\begin{lem}
Let $\lef x_n \rig \in l_p^w(X)$. Then $E_x (Ball(l_{p'}))$ is $p$-limited if and only if for every $\lef f_n \rig \in l_p^{w^*}(X^*)$, we have $\lef \lef f_n(x_k) \rig_k \rig_n \in l_p^s(l_p)$.
\end{lem}
\pf
Consider $x=\lef x_n\rig \in l_p^w (X)$ such that $S=E_x(Ball (l_{p'}))$ is $p$- limited. Then given $\lef f_n \rig \in l_p^{w^*}(X^*)$, there is an $\lef \alpha _n\rig \in l_p$ such that for each $\beta = \lef \beta _k \rig \in Ball (l_{p'})$ we have 
$$|f_n (E_x (\beta))|\leq \alpha _n, \ \textnormal{for all} \ n.$$
i.e.
$$|\lef \beta ,\lef f_n(x_k)\rig \rig _{k=1}^\infty | \leq \alpha _n, \ \textnormal {for all} \ n.$$
It follows that $\|\lef f_n (x_k) \rig _{k=1}^\infty \| _p \leq \alpha _n$, for all $n$. Thus $\lef\lef f_n (x_k) \rig _k \rig _n \in l_p^s(l_p)$ for all $\lef f_n \rig \in l_p^{w^*}(X^*)$. 

Now tracing back the proof, we can prove the converse also. \hfill $\triangle$
\vskip 8pt plus 0fill

Since $l_p^{w^*}(X^*)$ can be identified with $B (X, l_p)$, where to each $f \in \lef f_n \rig \in l_p^{w^*}(X^*)$ we get $(E_f)_* \in B(X, l_p)$ given by $(E_f)_* (x)= \lef f_n(x) \rig$ with $\|f\|_p^{w^*} = \|(E_f)_*\|$, and since $l_p^{w^*}(X^*) = l_p^w(X^*)$ the above lemma can be reorganized as follows.

\begin{pro}
Let $x = \lef x_n \rig$ be a sequence in $X$. The following are equivalent. \par
{\leftskip=6mm
\item{(a)} $x \in l_p^w(X)$ and $E_x(Ball(l_{p'}))$ is a $p$-limited set in $X$.
\item{(b)} $\lef Tx_n \rig \in l_p^s(l_p)$ for all $T \in B(X, l_p)$. 
\item{(c)} $E_x \in \Pi _p^d (l_{p^{\prime}}, X)$. \par}
\end{pro}

In this way we observe a new notion of summability in Banach spaces in the `$p$-sense'. Let us rename this phenomena as follows:

\begin{Def}
A sequence $\lef x_n \rig$ in $X$ is said to be operator $p$-summable in $X$ if it satisfies one (and hence all) of the conditions of Proposition 2.4.
\end{Def}

Note that every norm $p$- summable sequence in $X$ is operator $p$-summable. To see this, let $\lef x_n \rig \in l_p^s(X)$ and $T \in B(X, l_p)$.
Then $\|Tx_n\| \leq \|T\| \|x_n\|$ for all $n$ so that $\lef Tx_n \rig \in l_p^s(l_p)$. Thus $\lef x_n \rig $ is operator $p$-summable. We have already seen that an operator $p$-summable sequence is weakly $p$-summable.

\section{Towards weak summability}

In this section we characterize Banach spaces with the property that every weakly $p$-summable sequence is operator $p$-summable and give some examples of such spaces. We shall call a Banach space with this property a {\it weak $p$-space}. A simple characterization of such spaces in terms of operator ideals is given below. Let $X$ and $Y$ be Banach spaces and let $1 \leq p < \infty$. Then an operator $T \in B(X, Y)$ is called absolutely $p$- summing if for every $\lef x_n \rig \in l_p^w(X)$, $\lef Tx_n \rig \in l_p^s(Y)$. The set of all absolutely $p$-summing operators in $B(X, Y)$ is denoted by $\Pi _p(X, Y)$.

\begin{pro}
Let $X$ be a Banach space and $1 \leq p < \infty$. Then $X$ is a weak $p$-space if and
only if $\Pi _p(X, l_p) = B(X, l_p)$.
\end{pro}
\pf
Let $T \in B(X, l_p)$ and $x = \lef x_n \rig \in l_p^w(X)$. Suppose $X$ is a weak $p$-space.
Then $\lef x_n \rig$ is
operator $p$-summable so that $\lef Tx_n \rig \in l_p^s(l_p)$. Thus $T \in \Pi _p (X, l_p)$.
Tracing back we can prove the converse.  \hfill $\triangle$
\vskip 8pt plus 0fill
Before we give some examples of weak $p$- spaces we shall further explore Banach spaces
that satisfy an operator ideal equation of the above type. Given  Banach spaces $X$ and $Y$,
let $W(X, Y)$ and $\nu (X, Y)$ denote the sets of weakly compact and completely continuous
operators from $X$ to $Y$ respectively. Recall that a Banach space $X$ is said to have the
Dunford-Pettis property (D. P. P., for short) if for any Banach space $Y$,
$W(X, Y) \subset \nu (X, Y)$.

In 1940, Dunford and Pettis \cite{DP} proved that every weakly compact operator
defined on a $L_1(\mu)$ space takes weakly compact sets to norm compact sets. In 1953,
Grothendieck \cite{G2} defined a Banach space $X$ to have the \DPP if weakly compact
operators defined on $X$ are completely continuous and proved that $C(K)$ spaces also have this property. This result
was also obtained independently in 1955 by Bartle, Dunford and Schwartz \cite{BDS}.
Brace \cite{B} and Grothendieck \cite{G2} gave some nice characterizations of the
Dunford-Pettis property. A detailed survey of \DPP can be found in \cite{D}. In this
section we propose to extend this property to a $p$-setting to meet our above mentioned
end. For this purpose we recall the following characterization of D.P.P. essentially
due to Grothendieck \cite{G2}.

\begin{thm}
Let $X$ be a Banach space then the following statements are equivalent. \par
{\leftskip =6mm
\item {(a)} $W(X, Y) \subset \nu (X, Y)$ for all Banach spaces $Y$.
\item {(b)} $W(X, c_0) \subset \nu (X, c_0)$.
\item {(c)} For $\lef x_n \rig \in c_0^w(X)$ and $\lef f_n \rig \in c_0^w(X^*)$,
$\lef \lef f_k (x_n) \rig_k \rig_n \in c_0^s(c_0)$.
\item {(d)} For $\lef x_n \rig \in c_0^w(X)$ and $\lef f_n \rig \in c_0^w(X^*)$,
$\lef f_n (x_n) \rig \in c_0$. \par}
\end{thm}
Picking up (c) as an end, we now propose the following definition.

\begin{Def}
Let $1 \leq p, q \leq \infty$. A Banach space $X$ is said to have the $(p, q)$-Dunford-
Pettis property ($(p, q)$-D.P.P., for short) if given $\lef x_n \rig \in l_q^w(X)$
and $\lef f_n \rig \in l_p^w(X^*)$, we have $\lef \lef f_k(x_n)\rig_k \rig_n \in l_q^s(l_p)$.
For $p$(or $q$)$= \infty$, $l_p$ (or $l_q$) is replaced by $c_0$. For all $p$, the
$(p, p)$- Dunford-Pettis property shall be called the $p$- Dunford-Pettis property.
\end{Def}

It is immediate from Theorem 3.2 that the $\infty$-D.P.P. is the classical Dunford-Pettis
property.  In what follows we shall extend the above characterization theorem to the
$(p, q)$-setting. Towards this end the notion of weak $p$-compactness studied
by the authors \cite{SK} (also see Castillo \cite{C1, C2}) fits smugly in the scheme.
Let $X$ and $Y$ be Banach spaces and let $1 \leq p < \infty$.
An operator $T \in B(X, Y)$ is said to be $p$-compact (weakly $p$-compact) if
$T(Ball (X))$ is relatively $p$-compact (respectively, relatively weakly $p$-compact).
Here $l_{p'}$ is replaced by $c_0$ if $p = 1$. Let $W_p(X, Y)$ denote the set of all
weakly $p$-compact operators and $K_p(X, Y)$ that of all $p$-compact operators. The
next result was obtained by the authors \cite{SK}.

\begin{thm}
Let $X$ and $Y$ be Banach spaces, $1 \leq p < \infty$ and $ T \in B( X, Y)$. Then the
following statements are equivalent:\par
{\leftskip =6mm
\item {(a)} $T$ is weakly $p$-compact.
\item {(b)}There are $y \in l_p^w(Y) \ and \ S_y \in B( R(y), X^*)$ such that
$T^* = S_y \cdot E_y^*$, where $R(y) = \overline {Range (E_y^*)} \subset l_p$. \par}
\end{thm}

The set $W_p (X, Y)$ of all weakly $p$-compact operators in $B (X, Y)$ is a Banach
operator ideal with the factorization norm $\omega _p$ defined as follows:
$$\omega _p(T) = inf \{ \|S_y\| \|E_y\|: T^* = S_y \cdot E_y^* \ \
\textnormal {as in theorem 3.4(b)} \}.$$

\noindent Let $(A, \alpha)$ be an operator ideal. For Banach spaces $X$ and $Y$ we put
$$A^d(X, Y) = \{ T \in B(X, Y): T^* \in A(Y^*, X^*) \}.$$

\noindent For an operator $T \in A^d(X, Y)$, we put $\alpha ^d(T) = \alpha (T^*)$. With these notations
$(A^d, \alpha ^d)$ is also a Banach operator ideal and is called the dual ideal of $(A, \alpha)$.

\begin{cor}
For Banach spaces $X$ and $Y$, $T \in W_p^d(X, Y)$ if and only if there are $f = \lef f_n \rig
\in l_p^w(X^*)$ and $S_f \in B(R(f), Y)$ such that $T = S_f \cdot (E_f)_*$.
Here, $R(f) = \{ \lef f_n(x) \rig : x \in X \} \subset l_p$ and $(E_f)_* = E_f^*|_X$.
\end{cor}

\vskip 8pt plus 0fill
We can now extend the classical characterization theorem for the \DPP to the $(p, q)$-setting.
Let $X$ and $Y$ be a pair of Banach
spaces and $T \in B(X, Y^*)$. Then we have $T = i_Y^* \cdot T^{**} \cdot i_X$. Indeed,
for $x \in X$ and $y \in Y$ we have $\lef i_Y^* \cdot T^{**} \cdot i_X (x), y \rig =
\lef Tx, y \rig$. Here $i_X: X \hookrightarrow X^{**}$ is the canonical embedding.

\begin{thm}
Let $1 \le p < \infty$, $1 \le q \le \infty$ and $X$ a Banach space. Then the following
statements are equivalent:\par
\par {\leftskip=6mm
\item{(a)} $X$ has the $(p, q)$-Dunford-Pettis property.
\item{(b)} $W_p(Y, X^*) \subset \Pi _q^d(Y, X^*)$ for every Banach space $Y$.
\item{(c)} $\Pi _q^d(l_{p\prime}, X^*) = B(l_{p\prime}, X^*)$. \par}
\end{thm}
\pf
It only remains to show that $(a)$ implies $(b)$, for $W_p(l_{p\prime}, X^*) =
B(l_{p\prime}, X^*)$. To this end, assume that $X$ has the $(p, q)$-DPP and
let $T \in W_p(Y, X^*)$. Then by Theorem 3.4, there are ${\bf f} =
\lef f_n \rig \in l_p^w(X^*)$ and $S_{\bf f} \in B(R({\bf f}), Y^*)$ such that $T^* =
S_{\bf f} \cdot E_{\bf f}^*$. First we show that
$T^*|_X \in \Pi _q(X, Y^*)$. To see this, let $\lef x_n \rig \in l_q^w(X)$. Then as
$\lef f_n \rig \in l_p^w(X^*)$, we have $\lef (E_{\bf f})^* \cdot i_X(x_n) \rig =
\lef \lef f_k(x_n) \rig _k \rig _n \in l_q^s(l_p)$. Thus $\lef T^* \cdot i_X(x_n) \rig
\in l_q^s(Y^*)$. In other words, $T^* \cdot i_X \in \Pi _q(X, Y^*)$. It follows from
Proposition 2.19 in \cite{DJT} that $i_X^*\cdot T^{**} \in \Pi _q^d(Y^{**}, X^*)$, so that
$T = i_X^* \cdot T^{**} \cdot i_Y \in \Pi _q^d(Y, X^*)$. This completes the proof.  \hfill $\triangle$
\vskip 8pt plus 0fill

\noindent {\bf Note} In \cite{SK} the authors have observed that absolutely $p$-summing operators may be
regarded as $p$-completely continuous operators as they take weakly $p$-compact sets
to $p$-compact sets. Thus the classical \DPP may be traced back provided we regard
absolutely $p$-summing operators as $p$-completely continuous operators.
\vskip 8pt plus 0fill
In view of Proposition 3.1 and Theorem 3.6, we have the following characterization for weak $p$-spaces.

\begin{thm}
Let $1 \leq p < \infty$. Then for a Banach space $X$ the following statements are equivalent.
\par {\leftskip=6mm
\item{(a)} $X$ is a weak $p$-space
\item{(b)} $X$ has the $p$-Dunford-Pettis property.
\item{(c)} $W_p(Y, X^*) \subset \Pi _p^d(Y, X^*)$ for every Banach space $Y$.
\item{(d)} $\Pi _p^d(l_{p^\prime}, X^*) = B(l_{p^\prime}, X^*)$.
\item{(c$^\prime$)} $W_p^d(X, Y^*) \subset \Pi _p(X, Y^*)$ for every Banach space $Y$.
\item{(d$^\prime$)} $\Pi _p (X, l_p) = B(X, l_p)$ \par}
\end{thm}

\pf
Note that $W_p^d(X, l_p) = B(X, l_p)$. Thus in the light of  Proposition 3.1 and Theorem 3.6,
it is enough to show that (c)$\Leftrightarrow$(c$^\prime$) and that (d$^\prime$)$\Rightarrow$(d).
First assume that $W_p(Y, X^*) \subset \Pi _p^d(Y, X^*)$ and let $T \in W_p^d(X, Y^*)$. Then
$T^* \cdot i_Y \in W_p(Y, X^*) \subset \Pi _p^d(Y, X^*)$. Thus $i_Y^* \cdot T^{**} \in
\Pi _p(X^{**}, Y^*)$ so that $T = i_Y^* \cdot T^{**} \cdot i_X \in \Pi _p(X, Y^*)$  \cite [2.4]{DJT}.
Therefore, $W_p^d(X, Y^*) \subset \Pi _p(X, Y^*)$.
Next, let $W_p^d(X, Y^*) \subset \Pi _p(X, Y^*)$. If $T \in W_p(Y, X^*)$, then by
Theorem 3.4, $T^* = S_f \cdot E_f^*$ for some $f \in l_p^w(X^*)$. Thus by Corollary 3.5,
$T^* \cdot i_X \in W_p^d(X, Y^*) \subset \Pi _p(X, Y^*)$. It follows from Proposition 2.19
in \cite{DJT}, that $T = i_X^* \cdot T^{**} \cdot i_Y \in \Pi _p^d(Y, X^*)$. Thus
$W_p(Y, X^*) \subset \Pi _p^d(Y, X^*)$ so that (c)$\Leftrightarrow$(c$^\prime$). Now, let
$1 < p < \infty$ and assume that $\Pi _p(X, l_p) = B(X, l_p)$. Let $T \in B(l_{p^\prime}, X)$.
Then $T^* \cdot i_X \in B(X, l_p) = \Pi _p(X, l_p)$. Thus $i_X^* \cdot T^{**} \in
\Pi _p^d(l_{p^\prime}, X^*)$. As $l_{p^\prime}$ is reflexive, we have $T = i_X^* \cdot T^{**}$
so that $\Pi _p^d(l_{p^\prime}, X^*) = B(l_{p^\prime}, X^*)$. Finally, suppose that
$\Pi _1(X, l_1) = B(X, l_1)$. Let $T \in B(c_0, X^*)$. Then $T^* \cdot i_X \in B(X, l_1)
= \Pi _1(X, l_1)$. Thus $i_X^* \cdot T^{**} \cdot i_{c_0} \in \Pi _1^d(c_0, X^*)$ so that
$\Pi _1^d(c_0, X^*) = B(c_0, X^*)$. Therefore, (d$^\prime$)$\Rightarrow$(d), which completes
the proof.  \hfill $\triangle$

Some more consequences of Theorem 3.6 are in order.

\begin{cor} If $1 \leq q \leq p < \infty$ and if $X$ has the $(p, q)$-Dunford-Pettis
property, then it has the $p$-Dunford-Pettis property. In particular, $X$ is a weak
$p$-space.
\end{cor}

\begin{cor}
If $X^*$ has the $p$-Dunford-Pettis property, then so does $X$. In other words if $X^*$
is a weak $p$-space then so is $X$.
\end{cor}

\noindent {\bf Remark}: The $p= \infty$ case of the above corollary; i.e.,
if $X^*$ has the classical \DPP then so does $X$, was proved by
Grothendieck \cite{G2}.
\vskip 8pt plus 0fill
\noindent {\bf Note} It is interesting to note that Diestel et al. \cite[p. 433]{DJT}
defined a Banach space $X$ to be a Hilbert-Schmidt space if every Hilbert space operator
that factors through $X$ is a Hilbert-Schmidt operator. They observed that a Banach
space $X$ is a Hilbert-Schmidt space if and only if, for each $\lef x_n \rig \in l_2^w(X)$
and $\lef f_n \rig \in l_2^w(X^*)$, $\lef f_n(x_n) \rig \in l_2$. Clearly, all Banach
spaces satisfying the $2$-\DPP are Hilbert-Schmidt spaces. However, the class of all Hilbert-Schmidt
operators is clearly much larger. They proposed that, the class of Hilbert-Schmidt
spaces could be studied as the class of Banach spaces satisfying the \DPP of
``level 2". However, we clearly see that Banach spaces satisfying the $2$-\DPP are
more appropriate than the class of Hilbert-Schmidt spaces as the former exactly mimics
the classical case geometrically as well as analytically.

\vskip 8pt plus 0fill
\noindent {\bf Examples} 1. Let $X$ be an $\mathcal L_\infty$-space. If $1 \leq
p \leq 2$, then $X$ has the $(p, 2)$-DPP and $2$ is sharp \cite{G1, LP}. If
$2 < p < q < \infty$, then $X$ has the $(p, q)$-DPP and $q$ is sharp, that is
to say that $X$ has the almost $p$-DPP for every $p > 2$ \cite {R, K}.

\noindent 2. In view of Theorem 3.6 above, every $\mathcal L_1$-space has the above properties. In
particular, $c_0$ and $l_1$ have the $2$-DPP, the almost $p$-DPP if
$p >2$ and also the $\infty$-DPP (= \DPP).
\vskip 8pt plus 0fill
It is interesting to note that
these are the only $\mathcal L_p$-spaces with any $(r, s)$-DPP,
$1 \leq r, s \leq \infty$.

\begin{thm}
Let $1 < p < \infty$. Then $l_p$ does not have the $r$-\DPP for any $r > 1$.
In other words, for $1 <p < \infty$, $l_p$ is not a weak $r$-space for any $r > 1$.
\end{thm}

\pf
We divide the proof in several parts.

\noindent {\it Case 1}. Let $r \geq max \{p, p'\}$. Let $\{e_n\}$ be the standard
unit vector basis of $l_p$ and $\{f_n \}$ that of {$l_{p'}$. Then $\lef e_n \rig \in
l_r^w(l_p)$ and $\lef f_n \rig \in l_r^w(l_{p'})$. Since $\lef \lef f_k(e_n)
\rig _k \rig _n = \lef \lef \delta _k^n \rig _k \rig _n \notin l_r^s(l_r)$,
where $\delta _k^n$ is the Kronecker delta, we conclude that $l_p$ does not have
the $r$-DPP if $r \geq max\{p, p'\}$.
\vskip 8pt plus 0fill
Before we proceed to the other cases, we need to prove the following lemma.

\begin{lem}
Let $1 \leq s \leq p'$, where $p'$ is the harmonic conjugate of $p$, $1<p< \infty$.
Find $t > s$ such that $\frac {1}{s} - \frac {1}{p^\prime} = \frac {1}{t}$. Then for any
$\lef \alpha _n \rig \in l_t$, $\lef \alpha _n e_n \rig \in l_s^w(l_p)$.
\end{lem}

\noindent {\it Proof of the Lemma}: If $\lef \beta _n \rig \in l_{p^\prime}$, then
\begin{eqnarray*}
(\sum _{n=1}^\infty{| \lef \beta , \alpha _n e_n \rig |^s)}^{1/s} &=&
(\sum _{n=1}^\infty{| \alpha _n \beta _n |^s)}^{1/s}\\
 &\leq & (\sum _{n=1}^\infty{|\alpha _n|^t})^{1/t}
(\sum _{n=1}^\infty {|\beta _n|^{p^\prime}})^{1/{^\prime}} \\
&<& \infty .
\end{eqnarray*}
\noindent Thus $\lef \alpha _n e_n \rig \in l_s^w(l_p)$. \hfill $\triangle$

\vskip 8pt plus 0fill
Now we consider the other cases of the theorem.
\vskip 8pt plus 0fill
\noindent {\it Case 2}. Let $1 < r < min\{p, p'\}$. Find $t_1, t_2 >1$ such that
$\frac{1}{t_1} = \frac {1}{r} - \frac {1}{p'}$ and $\frac {1}{t_2} = \frac {1}{r}
- \frac {1}{p}$. Then $\frac {1}{t_1} + \frac {1}{t_2} = \frac {2}{r} -1 <
\frac {1}{r}$. Thus we can find $\lef \alpha _n \rig \in l_{t_1}$ and $\lef
\beta _n \rig \in l_{t_2}$ such that $\lef \alpha _n \beta _n \rig \notin l_r$.
Now by the above lemma $\lef \alpha _n e_n \rig \in l_r^w(l_p)$ and $\lef \beta _n
f_n \rig l_r^w(l_{p^\prime})$. But
$$\lef \lef \lef \beta _k f_k, \alpha _n e_n \rig \rig _k \rig _n
\notin l_r^s(l_r).$$
Thus $l_p$ does not have the $r$-DPP if $1 \leq r < min\{p, p^\prime\}$.
\vskip 8pt plus 0fill
\noindent {\it Case 3}. Let $r$ lie between $p$ and $p^\prime$. Note that $l_p$ has the
$r$-DPP if and only if $l_{p^\prime}$ has the $r$-DPP. Thus without any loss of generality
we may assume that $p <r < p'$. Find $t >1$ such that $\frac{1}{t} = \frac {1}{r}
- \frac {1}{p'}$. Then $r <t$ so that we can find $\lef \alpha _n \rig \in l_t$
with $\lef \alpha _n \rig \notin l_r$. Then $\lef \alpha _n e_n \rig \in
l_r^w(l_p)$. Also $\lef f_n \rig \in l_r^w(l_{p^\prime})$. But
$$\lef \lef \lef f_k, \alpha _n e_n \rig \rig _k \rig _n \notin l_r^s(l_r).$$
Thus $l_p$ does not have the $r$-DPP if $r$ lies between $p$ and $p^\prime$.
\vskip 8pt plus 0fill
Finally, since $\Pi _p(l_p) \neq  B(l_p)$, we conclude that both $l_p$ and $l_{p^\prime}$
do not have the $p$- and $p^\prime$-DPP. This completes the proof.  \hfill $\triangle$

\begin {cor}
For $1 < p < \infty$, $l_p$ does not have
\par {\leftskip=6mm
\item{(a)} The $(r, s)$-\DPP if $1 < r, s < \infty$.
\item{(b)} The $(r, 1)$-\DPP if $1 < r < \infty$.
\item{(c)} The $(1, r)$-\DPP if $1 < r < \infty$. \par}
\end{cor}

\pf
(a) For $1 < s \leq r$, if $l_p$ has the $(r, s)$-DPP then it also has the $r$-DPP.
Thus if $1 < s \leq r < \infty$, then $l_p$ does not have the $(r, s)$-DPP.

Next, let $1 < r < s$. Find $\lef \alpha _n \rig \in l_s^w(l_p)$ such that $\lef
\alpha _n \rig \notin l_r^w(l_p)$. Find $\beta \in l_{p^\prime}$ such that $\sum _{n=1}^\infty
|\lef \beta , \alpha _n \rig |^r = \infty$. Putting $\beta _1 = \beta$ and
$\beta _n = 0$ for $n \geq 2$, $\lef \beta _n \rig \in l_r^w(l_{p^\prime})$. However,
$$\lef \lef \lef \beta _k, \alpha _n \rig \rig _k \rig _n \notin l_s^s(l_r).$$
Thus $l_p$ does not have the $(r, s)$-DPP for $1 < r, s < \infty$.

Now, both (b) and (c) can be obtained on the lines of (a).
\hfill $\triangle$
\vskip 8pt plus 0fill
\noindent {\bf Note} We have not been able to settle whether for $1 < p < \infty$,
$l_p$ has the $1$-DPP.

\section{Towards norm summability}

In this section we shall examine a condition that forces every operator $p$-summable sequence to become norm-$p$-summable. Let $X$ be a Banach space and $1\leq p < \infty$. If $x \in l_p^w(X)$ is such that $E_x \in \Pi _p(l_{p^\prime}, X)$, then it follows from Proposition 5.5(a) in \cite{SK} and by Proposition 2.4 that $x$ is an operator $p$-summable sequence in $X$. In the light of this observation, we propose to study an operator version of the operator $p$-summable sequences.

\begin{Def}
An operator $T \in{\mathcal B}(X,Y)$ is said to be  $p$-limited if $T(Ball X)$ is $p$-limited in $Y$ and $T$ is said to be sequentially $p$-limited if $\lef Tx_n \rig$ is operator $p$-summable for all $\lef x_n \rig \in l_p^w(X)$.
\end{Def}
It follows from Proposition 2.4 that a sequence ${\bf x} = \lef x_n \rig $ in $X$ is operator $p$-summable if and only if $E_x \in {\mathcal B}(l_{p^\prime} , X)$ is a $p$-limited operator  if and only if $E_x \in \Pi _p^d(l_{p^\prime} , X)$. Further, we have
\begin{pro}
Let $1 \leq p < \infty$. Every $p$-limited operator $T \in {\mathcal B}(X,Y)$ is sequentially $p$-limited.
\end{pro}
\pf
Let ${\bf x}=\lef x_n \rig \in l_p^w(X)$. We may assume that $\left\|\lef x_n \rig\right\|_p^w \leq 1$, so that
$E_x(Ball (l_{p^\prime})) \subset Ball(X)$. Since $T(Ball (X))$ is $p$-limited in $Y$. Thus $T(E_x(Ball (l_{p^\prime})))$ is also $p$-limited in $Y$. Now by Lemma 2.2, $\lef Tx_n \rig$ is operator $p$-summable in $Y$.
\hfill $\triangle$

The following result will be used to characterize sequentially $p$-limited operator.
\begin{lem}
Let $1 \leq p < \infty$ and let $\alpha \in l_p^w(l_p)$. Then $\alpha \in l_p^s(l_p)$ if and only if $E_{\alpha} \in \Pi _p(l_{p^\prime}, l_p) = \Pi _p^d(l_{p^\prime}, l_p)$. Here $l_{p^\prime} = c_0$ when $p=1$.
\end{lem}
\pf
When $1 < p < \infty$, this fact follows from Remark (v) after Proposition 5.5 in \cite{SK}. Thus we may assume that $p = 1$. Again in this case, it follows, from \cite[Proposition 5.5(a)]{SK}, that if $E_\alpha \in \Pi _1(c_0, l_1)$, then $\alpha \in l_1^s(l_1)$.

Conversely, let $\alpha \in l_1^s(l_1)$, $\alpha = \lef \alpha _n \rig = \lef \lef \alpha _n^k \rig_k \rig_n$ such that $\alpha _n = \lef \alpha _n^k \rig_k \in l_1$ for all $n$. Put $\tilde {\alpha}_k = \lef \alpha _n^k \rig_n$ for all $k$. Then
$$\sum _{k=1}^\infty{ \sum _{n=1}^\infty{| \alpha _n^k|}} = \sum _{n=1}^\infty{ \sum _{k=1}^\infty{| \alpha _n^k|}} = \|\alpha \|_1^s,$$
so that $\tilde {\alpha}_k \in l_1$ for all $k$ with $\tilde{\alpha} = \lef \tilde{\alpha}_k \rig \in l_1^s(l_1)$. If $\beta = \lef \beta _n \rig = \lef \lef \beta _n^m \rig_m \rig_n \in l_1^w(c_0)$. Then
\begin{eqnarray*}
\lef E_\alpha(\beta _n)\rig_n &=& \lef \sum _{m=1}^\infty {\beta _n^m \alpha _m} \rig_n \\ &=& \lef \lef\sum _{m=1}^\infty {\beta _n^m \alpha _m^k} \rig_k \rig_n \\ &=& \lef \lef( \beta _n, \tilde{\alpha} _k \rig_k\rig_n.
\end{eqnarray*}
Since
$$ \sum _{k=1}^\infty{\sum _{n=1}^\infty {|(\beta _n, \tilde{\alpha}_k)|}} \leq \|\beta \|_1^w \sum _{k=1}^\infty{\|\tilde{\alpha}_k\|} = \|\beta \|_1^w \| \alpha \|_1^s,$$
we get $\lef E_\alpha (\beta _n)\rig_n \in l_1^s(l_1)$ with $\| \lef E_\alpha (\beta _n)\rig_n \|_1^s \leq \|\beta \|_1^w \|\alpha \|_1^s$. Thus $E_\alpha \in \Pi _1(c_0, l_1)$ with $\pi _1 (E_\alpha ) \leq \|\alpha \|_1^s$. Since $\| \alpha \|_1^s \leq \pi _1(E_\alpha)$, we conclude $\pi _1(E_\alpha) = \| \alpha \|_1^s$. \hfill $\triangle$
\begin{thm}
Let $T \in{\mathcal B}(X,Y)$. For $1 \leq p < \infty$, the following are equivalent:
\par {\leftskip=6mm
\item{(1)} $T$ is sequentially $p$-limited.
\item{(2)} $TU \in \Pi _p^d(l_{p^\prime},Y)$ for all $U \in{\mathcal B}(l_{p^\prime},X)$.
\item{(3)} $ST \in \Pi _p(X,l_p)$ for all $S \in {\mathcal B}(Y,l_p)$.\par}
\end{thm}
\pf
That (1) is equivalent to (2) follows from Lemma 2.2.

Now let (1) hold. If $S \in B(Y, l_p)$ and if $x = \left\langle x_n \right\rangle \in l_p^w(X)$ so that $E_x \in B(l_{p^{\prime}})$, then $T E_x \in \Pi _p^d(l_{p^{\prime}}, X)$.  Thus by Lemma 2.5, it follows that $STE_x \in \Pi _p^d(l_{p^{\prime}}, l_p) = \Pi _p(l_{p^{\prime}}, l_p)$. In other words, $\left\langle STx_n \right\rangle \in l_p^s(l_p)$ so that $ST \in \Pi _p(X, l_p)$. Thus (3) also holds. Finally, we can trace back the proof to show that (3) implies (1). \hfill $\triangle$

\noindent {\bf Remarks} \newline
1. If an operator $T\in \Pi _p(X,Y)\cup \Pi _p^d(X,Y)$, then $T$ is sequentially $p$-limited. \newline
2. Every sequentially $p$-limited operator in ${\mathcal B}(X, l_p)$ is in $\Pi _p(X, l_p)$.

Now we prove the following sequential characterization of subspaces of $L_p(\mu )$ whose operator characterization was obtained by Kwapi\'{e}n \cite{K}.
\begin{thm}
 Let $1 \leq p < \infty$. For a Banach space $X$ the following are equivalent:
\par {\leftskip=6mm
\item{(1)} Every operator $p$-summable sequence in $X$ is norm $p$-summable.
\item{(2)} $\Pi _p^d(Y,X) \subset \Pi _p(Y,X)$, for every Banach space $Y$.
\item{(3)} $\Pi _p^d(l_{p^\prime},X) = \Pi _p(l_{p^\prime},X)$.
\item{(4)} {\textrm(Kwapi\'{e}n)} $X$ is a subspace of $L_p(\mu)$ for some Borel measure $\mu$.\par}
\end{thm}
\pf
The equivalence of (2) and (4) was proved by Kwapi\'{e}n \cite{K}.

Let (1) hold. Assume that $T \in \Pi _p^d(Y, X)$ for some Banach space $Y$. If $y = \left\langle y_n \right\rangle \in l_p^w(Y)$, then $E_y \in B(l_{p^{\prime}}, Y)$. Thus $TE_y \in \Pi _p^d(l_{p^{\prime}}, X)$. Now, by Lemma 4 and assumption (1), we get that $\left\langle T y_n \right\rangle \in l_p^{op}(X) = l_p^s(X)$. It follows that $T \in \Pi _p(Y, X)$ so that (2) holds.

Since $\Pi _p(l_{p^{\prime}}, X) \subset \Pi _p^d(l_{p^{\prime}}, X)$ follows from \cite{SK}, we may conclude that (2) implies (3).

Finally, assume that (3) holds. Let $x \in l_p^{op}(X)$. Then $E_x \in \Pi _p^d(l_{p^\prime},X) = \Pi _p(l_{p^\prime},X)$. Now that $x \in l_p^s(X)$ again follows from \cite{SK}. This completes the proof. \hfill $\triangle$

Recall that a Banach space $X$ is said to have the Gelfand-Phillips property if every limited set in $X$ is relatively compact. Recall further that every limited set in a Banach space is conditionally weakly compact \cite{BD}. We do not know about the `$p$-version' of this result, possibly due to the absence of a $p$-prototype of a Rosenthal's $l_1$-theorem. At the same time let us note that in a Banach space, in which any operator $p$-summable sequence is norm $p$-summable, a (relatively) weakly $p$-compact set is (relatively) $p$-compact if and only if it is $p$-limited. Thus the condition that every operator $p$-summable sequence in a Banach space is norm $p$-summable can be seen as a $p$-version of the Gelfand-Phillips property.

\vskip 8pt plus 0fill
\noindent{\bf An operator ideal}: Let $1 \leq p < \infty$. For a pair of Banach spaces $X$ and $Y$, consider the set $Lt_p^\cdot (X, Y)$ of all sequentially $p$-limited operators in $B(X, Y)$. For $T \in Lt_p^\cdot (X, Y)$, we define
$$lt_p (T):= sup \{ \pi _p(ST): S\in B(Y, l_p) and \|S\| \leq 1\}.$$
Then it is a routine to prove the following result.
\begin{pro}
For $1 \leq p < \infty$, $(Lt_p^\cdot , lt_p)$ is a normed operator ideal.
\end{pro}

\noindent{\bf Note} We have not been able to show whether in general $Lt_p^\cdot(X, Y)$ is $lt_p$-complete. We, however, adopted the following approach.

Let $X$ and $Y$ be Banach spaces and $T\in B(X, Y)$.  For any $1 \leq p < \infty$, we can
define $\varphi _{T}: B(Y, l_p) \to B(X, l_p)$ given by $\varphi _{T}(S) = ST$, for all
$S \in B(Y, l_p)$. Now it is easy to show that $T \mapsto \varphi _{T}$ is a linear
isometry from $B(X, Y)$ into $B( B(Y, l_p), B(X, l_p))$; $1 \leq p \leq \infty$. Moreover,
if $1 \leq p < \infty$, it follows from Proposition 4.3, that $T \in Lt_p^\cdot(X, Y)$
if and only if $\varphi _T(B(Y, l_p)) \subset \Pi _p(X, l_p)$. In this case for all
$T \in Lt_p^\cdot (X, Y)$, we have
$$\\lt_p(T) = \textnormal {The operator norm of}\  \varphi _T\ \textnormal{in} \ B( B(Y, l_p), \Pi _p(X, l_p)).$$
We write $Lt _p(X, Y)$ for the completion of $\{\varphi _T: T \in Lt_p^\cdot(X, Y)\}$
in $B( B(Y, l_p), \Pi _p(X, l_p))$ and denote the operator norm on $Lt_p(X, Y)$ again by
$lt_p(.)$. Thus, proposition 4.11 may be re investigated in the following manner.

\begin{pro}
$(Lt_p, lt_p)$ is a Banach operator ideal, $1 \leq p < \infty$.
\end{pro}

\noindent {\bf Remarks} 1. If $x \in l_p^w(X)$ is operator $p$-summable, then
$\|x\|_p^w = \|E_x\| \leq lt_p(E_x) := lt_p(x)$. If $x \in l_p^s(X)$, then
$lt_p(x) \leq \|x\|_p^s$.

\noindent 2. For $T \in \Pi_p(X, Y)$, $\|T \| \leq lt_p(T) \leq \pi(T)$. If
$T \in \Pi _p(X, l_p)$, then $lt_p(T) = \pi _p(T)$.
\vskip 8pt plus 0fill

\end{document}